\documentclass[12pt,a4paper]{amsart}
\usepackage{euscript,amsfonts,amssymb,amsmath,amscd}

\usepackage[T2A]{fontenc}

\usepackage[utf8]{inputenc}        

\usepackage[english,russian]{babel}

\usepackage{hyperref}

\addto\captionsrussian{}
\addto\captionsrussian{}

\input{epsf}
\newcommand{\op}{\operatorname}
\newcommand{\m}{\mathbb}

\theoremstyle{plain}
\newtheorem{theorem}{Theorem}[section]
\newtheorem{lemma}[theorem]{Lemma}

\theoremstyle{definition}
\newtheorem{definition}[theorem]{Definition}

\begin{document}

\title{
Chern numbers of manifolds with torus action.
}
\author{Andrey Kustarev}


\begin{abstract}
We show that every set of numbers that occurs as the set of Chern numbers of an almost complex manifold $M^n$, $n\geqslant 6$, may be realized as the set of Chern numbers of a connected almost complex manifold $M^n$ with an almost complex action of two-dimensional compact torus. 
\end{abstract}

\maketitle

\section{Introduction}

There is a well-known problem (posed by Hirzebruch) asking which combinations of Chern numbers of almost complex manifolds (or, equivalently, which complex cobordism classes) may be realized by a~smooth connected complex projective variety. As shown by Milnor, any complex cobordism class may be realized by a disjoint union of smooth projective varieties, so the connectivity condition in the formulation of Hirzebruch problem is crucial. As shown in (\cite{geiges}), any complex cobordism class may be realized by a~connected smooth almost complex manifold. 

Chern numbers of compact almost complex and symplectic manifolds with compact Lie group action have been an object of intensive research in recent years (\cite{sabgod}, \cite{sab}, \cite{liliu}, \cite{pelayotolman}, \cite{tolmanweitsman}). Any smooth manifold $M$ with an~almost complex action having only isolated fixed points necessarily has $c_n(M)\geqslant 0$; many other restrictions are shown to exist in this case. In this paper we show that if fixed points are allowed to be non-isolated, then all restrictions on Chern numbers vanish.

\begin{theorem}
\label{main}
Every complex cobordism class of dimension $n\geqslant 6$ may be realized by a connected almost complex manifold $M^n$ with an almost complex action of two-dimensional compact torus. 
\end{theorem}

Moreover, the set of fixed points of the action consists only of isolated fixed points and two-dimensional oriented closed surfaces. 

\begin{definition}
A non-singular complex projective variety $M = M^n$ will be called {\it good} if it is a product of one or several varieties belonging to following classes:
\begin{itemize}
{\item compact complex curves of genus $g>1$,}
{\item non-singular compact projective toric varieties $X^n$,}
{\item Milnor hypersurfaces $H_{i, j}$ with $4\leqslant i\leqslant j$,}
\end{itemize}
such that there is at most one complex curve in the product. 
\end{definition}


\begin{lemma}
\label{good}
Any complex cobordism class may be represented by a~disjoint union of good varieties.
\end{lemma}

We will obtain the proof of theorem \ref{main} by applying the construction of fiber connected sum (see lemma \ref{fibersum}, cf. \cite{geigesbook} for contact structures) after establishing lemma~\ref{good}.

\section{Proof of lemma \ref{good}}

Recall that the complex cobordism ring $\Omega^*_U$ is generated by smooth stably complex manifolds modulo the equivalence relation: $M_1\sim M_2$ if and only if there exists stably complex manifold $N$ such that ${\partial N = M_1 \cup (-M_2)}$ where $(-M_2)$ is $M_2$ with reversed orientation and restriction of stably complex structure of $N$ to $M_1$ and $M_2$ is equivalent to stably complex structures on $M_1$ and $M_2$ respectively. The operations of taking disjoint union and product of manifolds turn $\Omega^*_U$ into a graded commutative ring; stably complex manifold $M^n$ determines a class $[M^n]\in \Omega^{-n}_U$.

As shown by Milnor and Novikov (\cite{milnor}, \cite{novikov}), the ring $\Omega^*_U$ is isomorphic to the polynomial ring $\m Z[a_1, a_2, \ldots]$, where $\deg a_{i} = -2i$. In particular, the manifold is zero in $\Omega^*_U$ if and only if all of its Chern numbers are zero. Moreover, the set of stably complex manifolds $\{M^{2n}\}, n\geqslant 1$, may solve as the set of polynomial generators for $\Omega^*_U$ if and only if for any~$n$ we have $s_n([M^{2n}]) = \eta(n)$, where
$$
\eta(n) = p, \mbox{~~~if } n + 1 = p^k \mbox{ for some prime } p,
$$
$$
\eta(n) = 1 \mbox{~~~otherwise},
$$
and $s_n$ is a {\it Milnor number} which is a Chern number corresponding to Newton symmetric polynomial of dimension $n$ (e.g. $s_1 = c_1$, ${s_2 = c_1^2-2c_2}$, $s_3 = c_1^3 - 3c_1c_2 + 3c_3$, etc.).

We proceed with the proof of lemma \ref{good} by induction. Any complex cobordism class of dimension two is represented by a disjoint union of complex curves, since any two-dimensional almost complex manifold is automatically a complex curve.

Suppose that the statement of theorem \ref{good} is true for all complex cobordism classes of dimension $\leqslant 2n$. 

\begin{lemma}
\label{snzero}
Any complex cobordism class $[M]\in \Omega^{-2n}_U$ satisfying ${s_n([M]) = 0}$ may be represented by a good variety. 
\end{lemma}

{\it Proof.}
Since $s_n([M]) = 0$, the class $[M]$ is decomposable and we have $[M] = \sum\limits_{k=1}^K [N_k\times N'_k]$ where $N_k$ and $N'_k$ are good varieties. If for some $k$ both $N_k$ and $N_k'$ contain complex curves of genus $g$ and $g'$ respectively, we replace these curves with disjoint unions of $(g-1)$ and $(g'-1)$ copies of $\m CP^1$ correspondingly. 
$\Box$

\begin{lemma}
If the class $[M]\in \Omega^{-2n}_U$ may be represented by a good variety, the same is true for $(-[M])$.
\end{lemma}

{\it Proof.} Recall that blowing up a point on a compact complex variety X lowers its Milnor number $s_n([X])$ by the value $(n+(-1)^n)$. This implies that in every dimension $n>1$ there exist toric varieties with positive and negative Milnor numbers -- namely, the complex projective space $\m CP^n$ and the variety obtained from $\m CP^n$ by blowing up three different fixed points of torus action. 

This argument shows that there exists a toric variety $N$ such that $s_n(a[M]+b[N])=0$, where $a,b\in \m Z_{>0}$. By lemma \ref{snzero}, the class $(-a[M]-b[N])$ may be represented by a disjoint union of good varieties. Then the same is true for the class ${-[M] = (-a[M]-b[N]) + b[N] + (a-1)[M]}$.
$\Box$

\begin{lemma}
For every $n\geqslant 1$ there exists a good variety $G_n$ (of complex dimension $n$) with $s_n([G_n])=\pm\eta(n)$.
\end{lemma}

{\it Proof.}
We will use the following result on toric generators in complex cobordism.

\begin{theorem}
\label{wilfong}
\cite{wilfong}
For every odd $n$ and every even $n\leqslant 100001$ there exists a compact non-singular toric variety $G_n$ (of complex dimension $n$) satisfying $s_n([G_n]) = \pm\eta(n)$. 
\end{theorem}

Therefore, it is enough to show that for every even $n>100001$ the complex cobordism generator $a_n\in\Omega^{-2n}_U$ may be realized by a disjoint union of good varieties. 

Recall that we have $s_n(\m CP^n) = n + 1$ and ${s_n(H_{i, j}) = -{i+j \choose i}}$ for $2\leqslant i \leqslant j$. By Kummer theorem, the maximum degree $k$ for which $p^k$ divides ${i+j \choose i}$ is equal to the number of  carries when $i$ is added to $j$ in base $p$. We need to show that ${\op{gcd}\Bigl(\left.\{(n+1), {i+j \choose i}\} \,\right\vert\, i + j = n+1,\, 2\leqslant i \leqslant j\Bigr) = \eta(n)}$ for every even $n$ large enough. 

Suppose first that $n + 1 = p^k$ and $\eta(n)=p$ for some prime $p$. Using \ref{wilfong}, we may only consider the case $p>2$. If $k=1$, then ${s_n([\m CP^n]) = n+1 = p = \eta(n)}$, so we may also assume that $k>1$.
Adding up numbers $p^{k-1}$ and $p^k - p^{k-1}$ in base $p$ has only one carry and by Kummer theorem we see that ${p^k \choose p^{k-1}}$ is divisible by $p$ and not by $p^2$. So the statement of the lemma is true unless $p^{k-1}\leqslant 3$. But in this case $p=3$ and $k=2$, so $n=8$ and toric generator $G_8$ exists by \ref{wilfong}.

Next, let us assume that $(n + 1)$ is not a degree of a prime $p$. If $p$ does not divide $(n+1)$, then there is nothing to prove since $s_n([\m CP^n])=n+1$. The representation of $(n+1)$ in base $p$ has at least two non-zero digits. Let $l>0$ be the maximum number such that $p^l<(n+1)$. Then ${n+1\choose p^l}$ does not divide by $p$. Recalling that $p$ divides $(n+1)$, we see that if $(n+1)-p^l\leqslant 3$, then either $n+1 = 2^l+2$ or $n+1=3^l+3$, so $n$ is odd and we may again apply \ref{wilfong}. 
$\Box$

\section{Fiber connected sums of almost complex manifolds}

In this section we apply lemma \ref{good} to prove the main result, theorem~\ref{main}. 

\begin{lemma}
\label{tordim}
Any disjoint union of good varieties admits a faithful action of compact torus $T^{\op{min}(4,n-1)}$ with the set of fixed points being a~disjoint union of nonsingular curves and isolated fixed points. 
\end{lemma}

{\it Proof.}
Recall that Milnor hypersurface $H_{i,j}\subset \m CP^i\times\m CP^j$, $i\leqslant j$, is given by the equation $\{z_0 w_0+\ldots+z_i w_i=0\}$ and has the natural faithful action of compact torus $T^i$ given by the formula
$$
(t_1,\ldots,t_i)\cdot (z_0,w_0, \ldots, z_i,w_i) = (z_0, w_0, t_1 z_1, t_1^{-1} w_1 , \ldots, t_i z_i, t_i^{-1} w_i).
$$
The action is faithful and has only isolated fixed points. Furthermore, every complex projective toric variety $X^n$ has a natural faithful action of compact torus $T^n$ and any curve of genus $g>1$ has a trivial torus action.

Taking products of torus actions, we see that every good variety $X^n$ possesses a faithful holomorphic action of a torus of dimension equal to at least $\op{min}(4, n-1)$. $\Box$

Let $X_1=X_1^{2n}$ and $X_2=X_2^{2n}$ be two smooth manifolds with a~smooth faithful action of torus $T^k$. Any principal orbit of the action has an~equivariant tubilar neighborhood diffeomorphic to ${T^k\times D^{2n-k}}$, where $D^{2n-k}$ is an open unit ball in $\m R^{2n-k}$ centered at the origin. Let~us choose two principal orbits $O_1\subset X_1$ and $O_2\subset X_2$ equipped with equivariant tubilar neigborhoods $U_1$ and $U_2$ respectively.  
One can now construct a {\it fiber connected sum} $Y = X_1\#_{T^k} X_2$ by gluing together diffeomorphic parts $U_1-O_1$ and $U_2-O_2$ of open manifolds $X_1-O_1$ and $X_2-O_2$ by the diffeomorphism $(t, r, \varphi)\sim (t, \frac1r, \varphi)$, where $t\in T^k$, $(r,\varphi)\in D^{2n-k}$, $r>0$, $\varphi\in S^{2n-k-1}$. 

There are smooth embeddings $(X_1-U_1)\hookrightarrow Y$ and ${(X_2-U_2)\hookrightarrow Y}$. If~$X_1$ and $X_2$ are almost complex manifolds, there is a natural question whether we can extend a $T^k$-invariant almost complex structure from a disjoint union of manifolds ${(X_1-U_1)\cup(X_2-U_2)}$ to the entire manifold $Y$.

\begin{lemma}
\label{fibersum}
Suppose that $X_1=X_1^{2n}$ and $X_2 = X_2^{2n}$ are compact manifolds with a smooth faithful almost complex action of compact torus $T^k$. If the stable homotopy group $\pi_{2n-k-1}(SO/U)$ is trivial, then the invariant almost complex structure can be extended from a disjoint union ${(X_1-U_1)\cup(X_2-U_2)}$ to an invariant almost complex structure on $Y$. Moreover, $[Y] = [X_1]+[X_2]$ in $\Omega^*_U$.
\end{lemma}

{\it Proof.}
The space $U_1-O_1\simeq U_2-O_2$ is a trivial fibration over the torus $T^k$. We denote by $A_t$ its fiber (an open annulus) over the point $t\in T^k$. Let us fix some point $t_0\in T^k$. Suppose that we managed to extend somehow the almost complex structure operator~$J$ from ${(X_1-U_1)\cup (X_2-U_2)}$ to the restriction of the tangent bundle $\tau_{\m R}(Y)|_{A_{t_0}}$. Then this new structure may be automatically extended by the action of torus $T^k$ to an invariant structure on the space ${U_1-O_1\simeq U_2-O_2}$, and therefore, on the entire space $Y$.

The obstruction to extending $J$ to $\tau_{\m R}(Y)|_{A_{t_0}}$ lies in a homotopy group $\pi_{2n-k-1}(SO(2n)/U(n))$. Since $k>0$, this group is equal to the stable homotopy group $\pi_{2n-k-1}(SO/U)$.

The equality $[Y]=[X_1]+[X_2]$ now follows from the localization theorem (\cite{ab}, \cite{toricgenera}, \cite{krichever}), which implies that Chern numbers of a manifold with an almost complex torus action are determined by the behaviour near fixed points.
$\Box$

{\it Proof of the theorem \ref{main}.}
By Bott periodicity, we have ${\pi_{j-1}(SO/U) = \pi_{j}(O)}$ and the stable homotopy group $\pi_j(O)$ is trivial if and only if $j \equiv 2,4,5,6 \op{mod} 8$. 

Consider the class $[M^n]\in \Omega^{-2n}_U$, $n\geqslant 3$. By lemma \ref{good} it can be realized by a disjoint union of good varieties and by lemma \ref{tordim} every of these varieties admits an action of torus $T^l$, where $l\geqslant \op{min}(4, n-1)$. Let $k = 4$, if $n \equiv 1(\op{mod} 4)$, and $k=2$ otherwise. We see that
\begin{itemize}
{\item $k\leqslant l$ and therefore every good variety admits a faithful almost complex action of torus $T^k$;}
{\item the stable homotopy group $\pi_{2n-k-1}(SO/U) = \pi_{2n-k}(O)$ is trivial.}
\end{itemize}
Applying lemma \ref{fibersum} finishes the proof of \ref{main}.
$\Box$

\end{document}